\documentclass[a4paper,12pt]{amsart}

\usepackage{amssymb,amsbsy,amsmath,amsfonts,amssymb,amscd}
\usepackage{latexsym}
\usepackage{amsxtra}
\usepackage{amscd}
\usepackage{graphics}
\usepackage{epic}
\usepackage{color}
\input xy
\xyoption{all}

\usepackage{ mathrsfs, amsfonts}

\newcommand\sA{{\mathcal A}}

\newcommand\sB{{\mathcal B}}

\newcommand\la{\lambda}

\DeclareMathOperator{\Hom}{Hom}

\def\Bbb{\bf}

\newcommand{\CC}{\ensuremath{\mathbb{C}}}

\newcommand{\ZZ}{\ensuremath{\mathbb{Z}}}
\newcommand{\QQ}{\ensuremath{\mathbb{Q}}}

\newcommand{\hol}{\ensuremath{\mathcal{O}}}

\newcommand{\PP}{\ensuremath{\mathbb{P}}}

\newcommand{\ra}{\ensuremath{\rightarrow}}

\def\eea{\end{eqnarray*}}
\def\bea{\begin{eqnarray*}}

\DeclareMathOperator{\Id}{Id}

\newcommand{\Proof}{{\it Proof. }}

\newcommand\dual{\mathrel{\raise3pt\hbox{$\underline{\mathrm{\thinspace d
\thinspace}}$}}}
\newcommand\qe{\ifhmode\unskip\nobreak\fi\quad $\Box$}       % box for QED

\def\BOX{\hfill\lower.5\baselineskip\hbox{$\Box$}}
\newcommand\Z{\Bbb Z}

\newtheorem{theo}[equation]{Theorem}
\newtheorem{remarkk}[equation]{Remark}
\newenvironment{rem}{\begin{remarkk}\rm}{\end{remarkk}}

\newtheorem{defin}[equation]{Definition}

\newtheorem{prop}[equation]{Proposition}

\newtheorem{lemma}[equation]{Lemma}
\newtheorem{example}[equation]{Example}

\newtheorem{question}[equation]{Question}

\newcommand{\sR}{\ensuremath{\mathcal{R}}}

\newtheorem*{example*}{Example}

\newcommand{\Spec}{{\rm Spec\,}}

\begin{document}

\title[SubcanonicalRings]{ Subcanonical graded rings which are not Cohen-Macaulay }
\author{ Fabrizio Catanese}
\address {Fabrizio Catanese: Lehrstuhl Mathematik VIII\\
Mathematisches Institut der Universit\"at Bayreuth\\
NW II,  Universit\"atsstr. 30\\
95447 Bayreuth, Germany\\}

\email{Fabrizio.Catanese@uni-bayreuth.de}

\address {Jonathan Wahl: Department of Mathematics, University of North Carolina
Chapel Hill NC 27599-3250, USA}
\email{jmwahl@email.unc.edu}

\thanks{
The present work took place in the realm of the DFG
Forschergruppe 790 ``Classification of algebraic
surfaces and compact complex manifolds''.}

\date{\today}

\maketitle
{\em  This article is  dedicated to  Rob  Lazarsfeld 
on the occasion of his $60$-th birthday.}

\begin{abstract}
We answer a question by Jonathan Wahl, giving examples of regular surfaces 
(so that the canonical ring is Gorenstein) with the following properties:

1)  the canonical divisor $K_S \equiv r L $
is a positive multiple of an ample divisor $L$ 

2)  the graded ring  $\sR : = \sR (X, L) $
associated to $L$ is not Cohen-Macaulay.

In the appendix Wahl shows how these examples lead to the existence of Cohen-Macaulay singularities with $K_X$  $\QQ$-Cartier which are not  $\QQ$-Gorenstein,
since their index one cover is not Cohen-Macaulay.

\end{abstract}

\tableofcontents
\section{Introduction}

The situation that we shall consider in this paper is the following: $L$ is an ample divisor on a complex projective manifold $X$ of complex dimension $n$,
and we assume that $L$ is subcanonical, i.e., there exists an integer $h$ such that we have the linear equivalence
$K_X \equiv h L$, where $h \neq 0$.

There are then two cases: $ h < 0$, and $X$ is a Fano manifold, or $ h > 0$ and $X$ is a manifold with ample canonical divisor (in particular $X$ is of general type).

Assume that $X$ is a Fano manifold, and that $ - K_X = r L$, with $ r > 0$: then, by Kodaira vanishing

$$ H^j (m L) : = H^j ( \hol_X (mL)) = 0 \  \forall m \in \ZZ , \forall  1 \leq j \leq n-1.$$

For $m < 0$ this follows from   Kodaira vanishing (and holds for $ j \geq 1$), while for $m \geq 0$ Serre duality gives
$ h ^j (m L) = h^{n-j} ( K - mL) =  h^{n-j} ( (-r  - m)L) = 0$.

At the other extreme, if  $K_X$ is ample, and $ K_X \equiv r L$, (thus $ r > 0$) by the same argument we get vanishing outside of the interval
$$ 0 \leq m \leq r. $$ 

To $L$ we associate as usual the finitely generated graded $\CC$-algebra 
$$ \sR (X, L) : = \oplus_{m \geq 0} H^0 ( X, \hol_X (mL))$$

Therefore in the Fano case, the divisor $L$ is arithmetically Cohen-Macaulay (see \cite{hartshorne}) and the above graded ring is a Gorenstein ring.

 The question is whether also in the case where $K_X$ is ample one may hope for such a good property.
 
 The above graded ring is integral over the canonical ring $\sA : = \sR (X, K_X) $, which is  a Gorenstein ring if and only if
 we have {\em pluri-regularity}, i.e., vanishing 
$$  H^j ( \hol_X ) = 0 \  \forall \  1 \leq j \leq n-1.$$
 
Jonathan Wahl asked the following question (which makes only sense for $ n \geq 2$):

\begin{question} {\bf (J. Wahl)}
Are there examples of subcanonical pluriregular varieties $X$ such that the graded ring $\sR (X, L) $ is not Cohen-Macaulay ?
\end{question}

We shall show that the answer  is positive, also in the case of regular subcanonical surfaces with $K_X$ ample, 
where by the assumption we have the vanishing
  $$ H^1 (m L) = 0 \  \forall m \leq 0, \ {\rm or} \  r \leq m $$
 and  the question boils down to
requiring the vanishing also for $ 1 \leq m \leq r-1$.

The following theorem answers the question by J. Wahl:
\begin{theo}\label{main}
For each $r = n-3$, where $n \geq 7$ is  relatively prime to $30$, and for each 
$  m,  \ 1 \leq m \leq r-1$ there are Beauville type surfaces $S$  with  $q(S) =0$ ($q(S) : =  \dim H^1 (S, \hol_S) $)
 s.t. $ K_S = r L$, and $ \  H^1 (m L) \neq  0 \ $.

\end{theo}

We get therefore examples of the following situation: $\sA : =  \sR (S, K_S)$ is a Gorenstein graded ring, and a subring of the ring $\sR : = \sR (S, L)$, which 
is not  arithmetically Cohen-Maculay; hence we have constructed 
 examples of non Cohen-Macaulay singularities ($\Spec (\sR )$) with $K_Y$  Cartier
which are cyclic quasi-\'etale covers of a Gorenstein singularity  ($\Spec (\sA )$).

In the Appendix, J. Wahl uses these to construct Cohen-Macaulay singularities with $K_X$  $\QQ$-Cartier whose index one cover is not Cohen-Macaulay.

In fact,  we can consider  three graded rings, two of which are subrings of the third, and which are cones 
associated to  line bundles on the surface S: 
\begin{itemize}
\item
 $Y : = \Spec (\sR )$, the cone associated to $L$, which is not  Cohen-Macaulay, while  $K_Y$ is  Cartier;
 \item
 $Z: = \Spec (\sA )$, the cone associated to $K_S$, which is Gorenstein;
 \item
 $X : = \Spec (\sB )$, the cone associated to $K_S + L$ (for instance), which is Cohen-Macaulay with $K_X$ $\QQ$--Cartier, but whose index 1 (or canonical) 
 cover $Y = \Spec (\sR )$ is not Cohen-Macaulay.
 \end{itemize}

\section{The special case of even surfaces}
Recall: a smooth projective surface $S$ is said to be  {\bf even } if there is a divisor $L$ such that
 $ K_S \equiv  2 L.$

This is a topological condition, it means 
 that 
the second Stiefel Whitney class $w_2(S) = 0$, or, equivalently, the 
intersection form 
$$  H^2(S, \ZZ) \ra \ZZ $$ is even (takes only even values).

In particular, an even surface is a minimal surface. 

In particular, if $S$ is of general type and even, the self intersection 

$$K^2_S = 4 L^2  = 8k$$
for some integer $k \geq 1$.

The first numerical case is therefore the case $K^2_S = 8$.

\begin{prop}\label{pg0}
Assume that $S$ is an even surface of general type with  $K^2_S = 8$  and $p_g(S) = h^0(K_S)=0$.
Then, if  $ K_S \equiv  2 L,$ then $ H^1 (L) = 0$.
\end{prop}

\Proof
We assume  that S is even, $K \equiv  2L$  , and $p_g = 0.$

Since  the intersection form is even, and $ K^2 \leq  9$ by the Bogomolov- Miyaoka - Yau inequality,  we obtain that $ L^2 = 2.$

The Riemann Roch theorem tells us:   $\chi(L) = 1 + \frac{1}{2} L (L-K) = 1 + \frac{1}{2} L (- L) = 0.$

On the other hand, by Serre duality $\chi(L) = 2 h^0(L) - h^1(L)$ , so if $H^1(L)$ is different from zero, then $ H^0(L) \neq 0$, contradicting $p_g = 0.$

\qed

Our construction for $n=5$ shall show in particular that the
`Beauville surface', constructed by  Beauville in \cite{beauville}  is an even surface with $K^2_S = 8, q (S) = p_g (S)= 0$,
but  with $  H^1 ( L) =  0 $.

\section{Canonical linearization on Fermat curves}

Fix a positive integer $ n \geq 5$, and let $C$ be the  degree $n$ Fermat curve
$$ C : = \{ (x,y,z) \in \PP^2 | f (x,y,z) : = x^n + y^n + z^n = 0 \}. $$ 

Let as usual $\mu_n$ be the group of $n$-roots of unity.

The the group $$G : = \mu_n^2 = \mu_n^3 / \mu_n$$ acts on $C$, and we obtain a natural  linearization
of $\hol_C(1)$ by letting  $(\zeta, \eta) \in \mu_n^2 $ act as follows:
$$ z \mapsto z , x   \mapsto  \zeta x , y   \mapsto  \eta y.$$  

In other words, $H^0 ( \hol_C(1))$ splits as a direct sum of one dimensional eigenspaces (respectively generated by $x,y,z$)
corresponding to the characters $ (1,0), (0,1), (0,0) \in (\ZZ/n)^2 \cong \Hom (G, \CC^*).$

Similarly, for $ m \leq n-1$,  the monomial $ x^a y ^b z^{m-a-b} \in H^0 ( \hol_C(m))$ generates the unique eigenspace for the character $(a,b)$
(we identify here $ \ZZ/n \cong \{0,1, \dots, n-1\}$ and we obviously require $ a+b \leq m$).

However, any two linearizations differ (see \cite{mumford}) by a character of the group. 

\begin{defin}
Assume that $n$ is not divisible by $3$.

We call the {\bf canonical } linearization on $H^0 ( \hol_C(1))$ the one obtained from the natural one by twisting with the character
$  (n-3)^{-1} (1,1)$. Thus $x$ corresponds to the character $ v_1 : = (1,0) + (n-3)^{-1} (1,1) = (-3)^{-1} (-2,1)$,
$y$ corresponds to the character $ v_2 : = (0,1) + (n-3)^{-1} (1,1) = (-3)^{-1} (1,-2)$,
$z$ corresponds to the character $ v_3 : = (-3)^{-1} (1,1) $.

\end{defin}

\begin{rem}

(I) Observe that $v_1, v_2$ are a basis of $(\ZZ/n)^2$ as soon as $n$ is not divisible by $3$.

Indeed, $v_1 + v_2 = \frac{1}{3} (1,1) = 3^{-1} (1,1)$, hence
$$  (1,0) = v_1 + 3^{-1}  (1,1) = 2 v_1 + v_2, \ \  (0,1) = 2 v_2 + v_1.$$ 

(II) Observe that the above linearization induces a linearization on all multiples of $L$,
and, in the case where $ m = (n-3)$, we obtain the natural linearization on the
canonical divisor of $C$, $ \hol_C (n-3) \cong \Omega^1_C$.

Since, if we take affine coordinates where $ z= 1$, and we let $ f$ the equation of $C$,
we have
$$ H^0 (\Omega^1_C) = \{ P(x,y) \frac{dx}{f_y}  =  - P(x,y) \frac{dy}{f_x} \}$$ 
and the monomial $P=  x^a y^b$ corresponds under this isomorphism to the character $ (a+1, b+1)$.

(III) In particular, Serre duality $$ H^0 (\hol_C(m)) \times H^1 (\Omega^1_C (-m)) \ra H^1 (\Omega^1_C ) \cong \CC,$$
where $\CC$ is the trivial $G$-representation, is $G$-invariant.
\end{rem}

From the previous discussion follows also

\begin{lemma}
The monomial  $ x^a y^b z ^c \in H^0( \hol_{C } (m))$ (here $ a,b,c \geq 0, \ a+ b + c = m$) corresponds to the character $\chi$ equal to:
$$(a,b) + (-3)^{-1} (m,m) =     (a - c   ) v_1 + ( b - c ) v_2 .$$
\end{lemma}

\Proof
$ v_1 + v_2 = \frac{1}{3} (1,1)$, hence $  (a,b) + (-3)^{-1} (m,m) = a v_1 + b v_2 + (-m+a+b) (3)^{-1}(1,1) =   (a - c   ) v_1 + ( b - c ) v_2 .$

\qed

\section{Abelian Beauville  Surfaces and their subcanonical divisors}

We recall now the construction ( see also \cite{isogenous}, or \cite{BCG}) of a Beauville surface with Abelian group $G \cong (\ZZ/n)^2$, 
where $n$ is not divisible by $2$ and by $3$.

\begin{defin}
(1) Let  $\Sigma \subset G$ be the union of the three respective subgroups generated by $(1,0), (0,1), (1,1)$.

(2) Let $\psi : G \ra G$ a homomorphism such that, setting $\Sigma^* : = \Sigma \setminus \{(0,0)\}$, $
\psi (\Sigma^*) \cap \Sigma^* = \emptyset$ (equivalently, $\psi (\Sigma) \cap \Sigma = \{(0,0)\}$).

(3) Let $C$ be the degree $n$ Fermat curve and let $$ S = (C\times C) / (\Id \times \psi)(G),$$
i.e., the quotient of $C \times C$ by the action of $G$ such that 
$ g ( P_1, P_2) = (g (P_1), \psi(g)(P_2))$.

\end{defin}

\begin{rem}
(i) By property (2) $G$ acts freely and $S$ is a projective smooth surface with ample canonical divisor.

(ii) The line bundle $\hol_{C \times C} (1,1)$ is $G \times G$ linearized, in particular it is $G\cong (\Id \times \psi)(G)$-linearized,
therefore it descends to $S$, and we get a divisor $L$ on $S$ such that the pull back of $\hol_S(L)$ is the
above $G$-linearized bundle. 

(iii) By the previous remarks, we have a linear equivalence $$ K_S \equiv (n-3) L. $$ 
\end{rem} 

\section{Cohomology of multiples of the subcanonical divisor $L$}

We consider now an integer $ m$ with
$$ 1 \leq m \leq n-4, $$

and we shall determine the space $ H^1 ( \hol_S (m L)) $. 

Observe first of all that  $ H^1 ( \hol_S (m L))  \cong  H^1( \hol_{C \times C} (m,m))^G$.

By the K\"unneth formula $$H^1( \hol_{C \times C} (m,m)) \cong $$
$$ \cong [H^0( \hol_{C } (m)) \otimes H^1( \hol_{C } (m)) ] \bigoplus 
 [H^1( \hol_{C } (m)) \otimes H^0( \hol_{C } (m)) ] .$$
 
 We want to decompose the right hand side as a representation of $G \cong  (\Id \times \psi)(G)$.
 
 Explicitly, $H^0( \hol_{C } (m)) = \oplus_{\chi} V_{\chi}$, where if we write the character $\chi = (a,b) + (-3)^{-1} (m,m)$
 ( $\chi  =    (a -  (m-a-b)  ) v_1 + ( b - (m-a-b) ) v_2 $ as we saw)
 then $V_{\chi}$ has dimension equal to one and corresponds to the monomial $ x^a y^b z ^{m-a-b}$,
 where $ a,b \geq 0, \ a+ b \leq m$.
 
 By Serre duality,  $H^1( \hol_{C } (m)) = \oplus_{\chi'} V_{- \chi'}$, where if we write as above $\chi ' = (a',b') + (-3)^{-1} (m',m')$,
 then $ V_{- \chi'}$ is the dual of  $ V_{ \chi'}$, corresponding to the monomial $ x^{a'} y^{b'} z ^{m'-a'-b'}$,
 where $ m' = n - 3 -m$, so $ 1 \leq m' \leq n-4 $ also, and where $ a',b' \geq 0, \ a'+ b' \leq m'$.

\bigskip

Now, the homomorphism $\psi : G \ra G$ induces a dual homomorphism $\phi : = \psi^{\vee} : G^{\vee} \ra G^{\vee} $,
therefore we can finally write $H^1( \hol_{C \times C} (m,m))$ as a representation of $G \cong  (\Id \times \psi)(G)$:
$$ H^1( \hol_{C \times C} (m,m)) = \bigoplus_{\chi, \chi'} [  ( V_{\chi} \otimes  V_{- \phi( \chi')}) \oplus (V_{- \chi'} \otimes V_{\phi( \chi)})].$$ 

We have proven therefore the

\begin{lemma}
$ H^1 ( \hol_S (m L)) \neq 0 $ if and only if there are characters $\chi = (a-c)  v_1 + (b - c)  v_2 $ and $\chi' = (a' - c') v_1 + (b' -c') v_2$
with $ a,b \geq 0, \ a+ b \leq m$, $ a',b' \geq 0, \ a'+ b' \leq m' = n-3-m$ such that 
$$ \chi = \phi (\chi') \ {\rm or} \   \chi' = \phi (\chi) .$$ 

\end{lemma}

\bigskip

{\em Proof of  theorem  \ref{main}.}

We take now $\phi$ to be given by a diagonal matrix in the basis $v_1, v_2$, i.e., such that 
$$ \phi (v_j) = \la_j v_j , \ j = 1,2,  \  \la_j \in (\ZZ/n)^*.$$

For further use we also set  $ \la : = \la_1, \mu : = \la_2. $ 

\bigskip

Given $n$ relatively prime to $30$  and $ 1 \leq m \leq n-4$, we want to find $\la_1$ and $ \mu$ such that the equations

$$ (a-c) = \la  (a' - c')  $$ 
$$ (b-c) = \mu (b' - c')  $$ 

have solutions with $a,b,c \geq 0$, $ a + b + c = m$, and  $a',b',c' \geq 0$, $ a' + b' + c' = m'$.

\bigskip

The first  idea is simply to take $b=c$ and $b' = c'$, so that $\mu$ can be taken arbitrarily.

For the first equation some care is needed, since we want that $\la$ be a unit: for this it suffices that 
$ (a-c) ,   (a' - c')  $ are both units, for instance they could be chosen to be equal  to one of the three numbers $1, 2,3$, according to the congruence class
of $m$, respectively $m'$, modulo $3$.

With this proviso we have to verify that we have a free action on the product.

\begin{lemma}
If $ n \geq 7$,  given $\la$ a unit, there exists a unit $\mu$ such that $\psi = \phi^{\vee} $ satisfies the condition $\psi (\Sigma) \cap \Sigma = \{(0,0)\}$.
\end{lemma}

\Proof
Since  $ (1,0) = 2 v_1 +  v_2$ and   $ (0,1) =  v_1 +  2 v_2$, the matrix of $\phi$ in the standard basis is the matrix 
 
 \[ \phi  = \frac{1}{3}  
\begin{pmatrix}
    4 \la - \mu&    2 (\la - \mu )\\
  2 (\mu - \la)   & 4 \mu -  \la
  \end{pmatrix} 
\] 

while the matrix of $\psi$ is the matrix 

 \[ \psi  = \frac{1}{3}  
\begin{pmatrix}
  A : =  4 \la -  \mu&   B : =  2 (\mu - \la)  \\
 C : = 2 (\la - \mu )   & D : =  4 \mu -  \la
  \end{pmatrix} 
\] 

The conditions for a free action boil down to:

$$ A, B, C, D,  A + B, C + D $$ are units in $\ZZ/ n$, 
and moreover $ A \neq B$, $C \neq D, A + B \neq  C + D$.

These are in turn equivalent to the condition that 

$$  \la , \mu , \la - 4 \mu , \la -  \mu, \mu - 4 \la, \la + 2 \mu , 2 \la + \mu \in (\ZZ/n)^*.$$ 

Given $\la \in (\ZZ/n)^*$, consider its direct sum decomposition given by the Chinese remainder theorem and the primary factorization of $n$. For each prime $p$ dividing $n$,
the residue classes modulo $p$ which are excluded by  the above condition are at most five values   inside $ (\ZZ/p)^*$,
hence we are done if $ (\ZZ/p)^*$ has at least six elements.

Now, since $n$ is relatively prime to $30$, each prime number dividing it is greater or equal to  $p=7$.

\qed

\begin{prop}
Consider the Beauville surface $S$ constructed in \cite{beauville}, corresponding to the case $n=5$.

Then $S$ is an even surface and $ K_S \equiv 2L$, where $ H^1 (L) = 0$.

\end{prop}

\Proof
We observe that $L$ is unique, because the torsion group of $S$ is of exponent $5$ (see \cite{bc}).

The existence of $L$ follows exactly as in the proof of the main theorem, where the condition $n \geq 7$ was not used.
That $ H^1 (L) = 0$ follows directly from proposition \ref{pg0}.

\qed

{\bf Acknowledgements.}
I would like to thank Jonathan Wahl for asking the above question. In the appendix below he describes  a construction
based on our main result.

%%%%%%%%%%%%%%%%%%%%%%%%%%%%%%%%%%%%%%%%%%%%%%%%%%%%%%

\section{Appendix by Jonathan Wahl: A non-$\QQ$-Gorenstein Cohen-Macaulay cone $X$ with $K_X$ $\QQ$-Cartier}
A germ $(X,0)$ of an isolated normal complex singularity of dimension $n\geq 2$ is called \emph{$\QQ$-Gorenstein} if 
\begin{enumerate}
\item $(X,0)$ is Cohen-Macaulay
\item The dualizing sheaf $K_X$ is $\QQ$-Cartier (i.e., the invertible sheaf $\omega_{X-\{0\}}$ has finite order $r$)
\item The corresponding cyclic \emph{index one} (or \emph{canonica}l) cover $(Y,0)\rightarrow (X,0)$  is Cohen-Macaulay, hence Gorenstein.
\end{enumerate}
  Alternatively, $(X,0)$ is the quotient of a Gorenstein singularity by a cyclic group acting freely off the singular point.  Some early definitions did not require the third condition, which is of course automatic for $n=2$. 

 If $(X,0)$ is $\QQ$-Gorenstein, a one-parameter deformation $(\mathcal X,0)\rightarrow (\CC,0)$ is called $\QQ$-Gorenstein if it is the quotient of a deformation of the index one cover of $(X,0)$; this is exactly the condition that $(\mathcal X,0)$ is itself $\QQ$-Gorenstein.  These notions were introduced by Koll\'ar and Shepherd-Barron \cite{ksb}, who made extensive use of the author's explicit smoothings of certain cyclic quotient surface singularities in \cite{lw} (5.9); these deformations were patently $\QQ$-Gorenstein, and it was important to name this property. 

 Recently, the author and others considered rational surface singularities admitting a rational homology disk smoothing (i.e., with Milnor number $0$). The three-dimensional total space of the smoothing had a rational singularity with $K$ $\QQ$-Cartier, but it was not initially clear whether the smoothings were $\QQ$-Gorenstein.  (This was later established \cite{wah}  by proving the stronger result that the total spaces were log-terminal.)  In fact,  one needs to be careful because of the examples of A. Singh:

\begin{example*} \cite{si}: There is a three-dimensional isolated rational (hence Cohen-Macaulay) complex singularity $(X,0)$ with $K_X$ $\QQ$-Cartier which however is not $\QQ$-Gorensteiin.
\end{example*}

The purpose of this note is to use F. Catanese's result to provide other examples; they are not rational, but are cones over a smooth projective variety,  which could for instance be assumed to be projectively normal with ideal generated by quadrics.

\begin{prop} Let $S$ be a surface as in Theorem 2 of Catanese's paper, with $h^1(S,\mathcal O_S)=0$, $L$ ample, $K_S=rL$ (some $r>1$), and $h^1(mL)\neq 0$ for some $m>0$.  Let $t$ be greater than $r$ and relatively prime to it.  Then \begin{enumerate}
\item The cone 
$R= \sR (S, tL) : = \oplus_{m \geq 0} H^0 ( S, \hol_S (mtL))$ is Cohen-Macaulay.
\item The dualizing sheaf of $ R$ is torsion, of order $t$.
\item The index one cover is  $ \sR (S, L) : = \oplus_{m \geq 0} H^0 ( S, \hol_S (mL))$, and is not Cohen-Macaulay.

\end{enumerate}
In particular, $R$ is not $\QQ$-Gorenstein.
\begin{proof} The Cohen-Macaulayness for $R$ follows because $h^1(itL)=0$, all $i$, thanks to Kodaira Vanishing.  Let $\pi:V\rightarrow S$ be the geometric line bundle corresponding to $-tL$; then
 $H^0(V,\mathcal O_V)\equiv R$.  Since $K_V\equiv \pi^*(K_S+tL)$, one has that $jK_R\equiv \oplus_{n\in \Z}H^0(S,j(K_S+tL)+nL)$; since $tK_S=r(tL)$ with $r$ and $t$ relatively prime, $K_R$ has order $t$.  Making a cyclic $t$-fold cover and normalizing gives that $\mathcal R(S,L)$ is the index one cover, which as Catanese has noted is not Cohen-Macaulay.
\end{proof}
  \end{prop}

\end{document}